\documentclass[12pt]{amsart}
\usepackage{amscd,amssymb}
\usepackage[graph,frame,poly,arc]{xy}  
\usepackage[plainpages,backref,urlcolor=blue]{hyperref}

\theoremstyle{plain}
\newtheorem{thm}[subsection]{Theorem}
\newtheorem{lem}[subsection]{Lemma}
\newtheorem{prop}[subsection]{Proposition}
\newtheorem{cor}[subsection]{Corollary}

\theoremstyle{definition}
\newtheorem{rk}[subsection]{Remark}
\newtheorem{definition}[subsection]{Definition}
\newtheorem{ex}[subsection]{Example}

\numberwithin{equation}{section}
\newcommand{\HH}{{\mathcal H}}

\newcommand{\A}{{\mathcal A}}
\newcommand{\B}{{\mathcal B}}

\newcommand{\D}{{\mathcal D}}
\newcommand{\CC}{{\mathcal C}}

\newcommand{\PPP}{{\mathcal P}}

\newcommand{\C}{\mathbb{C}}
\newcommand{\PP}{\mathbb{P}}

\newcommand{\N}{\mathbb{N}}

\DeclareMathOperator{\rank}{rank}

\DeclareMathOperator{\mult}{mult}

\begin{document}

\title [Curve arrangements, pencils, and Jacobian syzygies]
{Curve arrangements, pencils, and Jacobian syzygies}

\author[Alexandru Dimca]{Alexandru Dimca$^1$}
\address{Univ. Nice Sophia Antipolis, CNRS,  LJAD, UMR 7351, 06100 Nice, France. }
\email{dimca@unice.fr}

\thanks{$^1$ Partially supported by Institut Universitaire de France.}



\begin{abstract} 
Let $\mathcal C :f=0$ be a curve arrangement in the complex projective plane. If $\mathcal C$ contains a curve subarrangement consisting of at least three members in a pencil, then one obtains an explicit syzygy among the partial derivatives of the homogeneous polynomial $f$. In many cases this observation reduces the question about the freeness or the nearly freeness of $\mathcal C$ to an easy computation of Tjurina numbers. Some consequences for Terao's conjecture in the case of line arrangements are also discussed as well as the asphericity of some complements of geometrically constructed free curves.

\end{abstract}
 
\maketitle


\section{Introduction} 

Let $S=\C[x,y,z]$ be the graded polynomial ring in the variables $x,y,z$ with complex coefficients and let $\CC:f=0$ be a reduced curve of degree $d$ in the complex projective plane $\PP^2$. The minimal degree of a Jacobian syzygy for $f$ is the integer $mdr(f)$
defined to be the smallest integer $r \geq 0$ such that there is a nontrivial relation
\begin{equation}
\label{rel_m}
 af_x+bf_y+cf_z=0
\end{equation}
among the partial derivatives $f_x, f_y$ and $f_z$ of $f$ with coefficients $a,b,c$ in $S_r$, the vector space of  homogeneous polynomials of degree $r$.  The knowledge of the invariant $mdr(f)$ allows one to decide if the curve $\CC$ is free or nearly free by a simple computation of the total Tjurina number $\tau(\CC)$, see \cite{duPCTC}, \cite{Dmax}, and Theorem \ref{thmTAN} and Theorem \ref{thmPEN} below for nice geometric applications.

When $\CC$ is a free (resp. nearly free) curve in the complex projective plane $\PP^2$, such that $\CC$ is not a  union of lines  passing through one point,  then the exponents of $\A$ denoted by $d_1 \leq d_2$ satisfy  $d_1 =mdr(f)\geq 1$  and
one has
\begin{equation}
\label{sum1}
d_1+d_2=d-1,
\end{equation}
(resp. $d_1+d_2=d$). Moreover, all the pairs $d_1,d_2$ satisfying these conditions may occur as exponents, see \cite{DStexpo}. For more on free hypersurfaces and free hyperplane arrangements see 
\cite{KS}, \cite{OT}, \cite{T}, \cite{HS}.
A useful result is the following.

\begin{thm}
\label{thmF} Let $\CC:f=0$ be a reduced plane curve of degree $d$.
If $ r_0 \leq mdr(f)$ for some integer $r_0\geq 1$, then 
$$\tau(\CC)\leq (d-1)^2-r_0(d-r_0-1)$$
and equality holds if and only if $\CC$ is free with exponents $(r_0,d-r_0-1)$. 
In particular, the set $F(d,\tau)$ of free curves in the variety $C(d,\tau)\subset \PP(S_d)$ of reduced plane curves of degree $d$
with a fixed global Tjurina number $\tau$ is a Zariski open subset.
\end{thm}
The interested reader may state and prove the completely similar result for nearly free curves.
If the curve $\CC$ is reducible, one calls it sometimes a curve arrangement. When the curve $\CC$ can be written as the union of at least three members of one pencil of curves, we say that $\CC$ is a curve arrangement of pencil type. Such arrangements play a key role in the theory of line arrangements, see for instance \cite{FY}, \cite{PS} and the references therein.

In this note we show that the existence of a subarrangement $\CC'$ in a curve arrangement $\CC:f=0$, with $\CC'$ of pencil type, gives rise to an explicit Jacobian syzygy for $f$.
We start with the simplest case, when $\CC=\A$ is a line arrangement and the pencil type subarrangement $\CC'$ comes from an intersection point  having a high multiplicity, say $m$, in $\A$. This case was considered from a different point of view in the paper \cite{FV} by D. Faenzi
and J. Vall\`es. However, the construction of interesting syzygies from points of high multiplicity in $\A$ is
very explicit and elementary in our note, see the formula \eqref{syz1}, while in  \cite{FV} the approach involves a good amount of Algebraic Geometry. This explicit construction allows us to draw some additional conclusions for the nearly free line arrangements as well.

Our first main result  is the following.

\begin{thm}
\label{thm1} If $\A:f=0$ is  a line arrangement and $m$ is the  multiplicity of one of its intersection points,  then either $mdr(f) =d-m$, or 
$mdr(f) \leq d-m-1$ and then one of the following  two cases occurs.

\begin{enumerate}

\item  $mdr(f) \leq m-1$. Then equality holds, i.e. $mdr(f) = m-1$, one has the inequality $2m < d+1$ and the line arrangement $\A$ is free with exponents $d_1=mdr(f)=m-1$ and $d_2=d-m$;

\item  $m \leq mdr(f) \leq d-m-1$, in particular $2m < d$.
\end{enumerate}

\end{thm}
Theorem \ref{thmF} can be used to identify the free curves in the case (2) above.
We show by examples in the third section  that all the cases listed in Theorem \ref{thm1} can actually occur inside the class of free line arrangements. A number of corollaries of Theorem \ref{thm1} on the maximal multiplicity $m(\A)$ of points in a free or nearly free line arrangement $\A$ are given in the second section.

On the other hand, as a special case of a result in \cite{DSa} recalled in subsection \ref{Saito} below, we have the following.

\begin{prop}
\label{propDSa} If $\A:f=0$ is  a line arrangement,  then the maximal multiplicity $m(\A)$ of points in $\A$ satisfies
$$m(\A) \geq \frac{2d}{mdr(f)+2}.$$
In particular, if $\A$ is  free or nearly free with exponents $d_1 \leq d_2$, then
$$m(\A) \geq \frac{2d}{d_1+2}.$$
 
\end{prop}

This inequality is sharp, i.e. an equality, for some arrangements, see Example \ref{ex5}.

\medskip

We say that Terao's Conjecture holds  for a free hyperplane arrangement $\A$ if any other hyperplane arrangement $\B$, having an  isomorphic intersection lattice  $L(\B)=L(\A)$,  is also free, see \cite{OT}, \cite{Yo}. This conjecture is open even in the case of line arrangements in the complex projective plane $\PP^2$, in spite of a lot of efforts, see for instance  \cite{A1},  \cite{A2}. For line arrangements, since the total Tjurina number $\tau(\A)$ is determined by the intersection lattice $L(\A)$, it remains to check that $\A:f=0$ and $\B:g=0$ satisfy $mdr(f)=mdr(g)$ and then apply \cite{duPCTC}, \cite{Dmax}. 
Theorem 3.1 in \cite{FV} and our results above imply the following fact, to be proved in section 3.

\begin{cor}
\label{cor2} Let $\A$ be a free line arrangement with exponents $d_1 \leq d_2$. If $$m:=m(\A) \geq d_1,$$
then Terao's Conjecture holds for the line arrangement $\A$. In particular, this is the case when one of the following conditions hold.

\begin{enumerate}

\item $d_1=d-m$;

\item  $m \geq d/2$;

\item $d_1 \leq \sqrt {2d+1} -1$.
\end{enumerate}
 
\end{cor}

\begin{rk}
\label{rkY} (i) The fact that the Terao's Conjecture holds for the line arrangement $\A$ when $m=m(\A) \geq d_1+2$ was established in \cite{FV2} by an approach not involving Jacobian syzygies. The result for $m=m(\A) \geq d_1$ is implicit in \cite{FV}, see Theorem 3.1 coupled with Remarks 3.2 and 3.3. Moreover, the case $m=m(\A)= d_1-1$ for some real line arrangements is discussed in Theorem 6.2 in \cite{FV}.

(ii) The cases $d_1=d-m$ and $m \geq d/2$ in Corollary \ref{cor2} follow also from the methods described in \cite{Yo}, see in particular 
Proposition 1.23 (i) and Theorem 1.39. Corollary \ref{cor2} follows also from \cite{A1},  Theorem 1.1., the claims (1) and (3). 

(iii) The case (3) in Corollary \ref{cor2} improves Corollary 2.5 in \cite{Dmax} saying that Terao's conjecture holds for $\A$ if $d_1 \leq \sqrt {d-1}$.

\end{rk}

It is known that Terao's Conjecture holds for the line arrangement $\A$ when $d=|\A|\leq 12$, see \cite{FV}.  This result and the case (3) in Corollary \ref{cor2} imply the following.

\begin{cor}
\label{cor3} Let $\A$ be a free line arrangement with exponents $d_1 \leq d_2$. If $$d_1 \leq 4,$$
then Terao's Conjecture holds for the line arrangement $\A$. 
 \end{cor}
The stronger result, corresponding to $d_1 \leq 5$, is established in \cite[Corollary 5.5]{A1}.
In the case of nearly free line arrangements we have the following  result, which can be proved by the interested reader using the analog of Theorem \ref{thmF} for nearly free arrangements.

\begin{cor}
\label{corNF} Let $\A$ be a nearly free line arrangement with exponents $d_1 \leq d_2$. If $$m(\A) \geq d_1,$$
then any other line arrangement $\B$, having an  isomorphic intersection lattice  $L(\B)=L(\A)$,  is also nearly free.

\end{cor}

Now we present our results for curve arrangements. First we assume that $\CC$ is itself an arrangement of pencil type.

\begin{thm}
\label{thm11}  Let $\CC:f=0$ be a curve arrangement in $\PP^2$ such that the defining equation has the form
$$f=q_1q_2 \cdots q_m,$$
for some $m \geq 3$, where $\deg q_1=\cdots=\deg q_m=k\geq 2$ and  the curves $\CC_i:q_i=0$ for $i=1,...,m$ are members of the pencil $\PPP: u\CC_1+v\CC_2$. 
Assume that $\PPP$ has a $0$-dimensional base locus and that it contains only reduced curves. 
Then either $mdr(f) =2k-2$, or $m=3$, $mdr(f) \leq 2k-3$ and in addition one of the following  two cases occurs.

\begin{enumerate}

\item   $k \geq 4$ and  $mdr(f) \leq k+1$. Then equality holds, i.e. $mdr(f) = k+1$,  and the curve arrangement $\CC$ is free with exponents $d_1=k+1$ and $d_2=2k-2$;

\item   $k \geq 5$ and $k+2\leq  mdr(f) \leq 2k-3$.

\end{enumerate}

\end{thm}

Using \cite{duPCTC}, \cite{Dmax} we get the following consequence.

\begin{cor}
\label{corPen1} Let $\CC$ be a curve arrangement of pencil type such that the corresponding  pencil $\PPP: u\CC_1+v\CC_2$  has a $0$-dimensional base locus and that it contains only reduced curves.
If the number $m$ of pencil members, which are curves of degree $k$, is at least $4$, then $mdr(f) =2k-2$. In particular, in this case the curve $\CC$ is free  if and only if
$$\tau(\CC)=(d-1)^2-2(k-1)(d-2k+1),$$
resp. $\CC$ is nearly free if and only if
$$\tau(\CC)=(d-1)^2-2(k-1)(d-2k+1)-1.$$
\end{cor}
And again Theorem \ref{thmF} can be used to identify the free curves in the case (2) above.
Now we discuss the case of a curve arrangement containing a subarrangement of pencil type.

\begin{thm}
\label{thm13}  Let $\CC:f=0$ be a curve arrangement in $\PP^2$ such that the defining equation has the form
$$f=q_1q_2 \cdots q_mh,$$
for some $m \geq 2$, where $\deg q_1=\cdots=\deg q_m=k\geq 1$ and the curves $\CC_i:q_i=0$ for $i=1,...,m$ are  reduced members of the pencil  $u\CC_1+v\CC_2$. Assume that the curves $\CC_1:q_1=0$, $\CC_2:q_2=0$ and $\HH:h=0$ have no intersection points and that the curve $\HH$ is irreducible.
Then either $mdr(f) =2k-2+\deg(h)=d-(m-2)k-2$, or  $mdr(f) \leq d-(m-2)k-3$  and then one of the following  two cases occurs.

\begin{enumerate}

\item   $mdr(f) \leq (m-2)k+1$. Then equality holds, i.e. $mdr(f) =(m-2)k+1 $,  and the curve arrangement $\CC$ is free with exponents $d_1=(m-2)k+1$ and $d_2=d-(m-2)k-2$;

\item   $ (m-2)k+2 \leq mdr(f) \leq d-(m-2)k-3$.

\end{enumerate}

\end{thm}
In fact, this result holds also when $k=1$ and $\HH$ is just reduced, see Remark \ref{rklines2}.

\begin{rk}
\label{rklines1}

Note that when $\CC$ is a line arrangement, containing strictly the pencil type arrangement $\CC'$ and such that $\deg h >1$, (i.e. $\CC$ contains at least two lines not in $\CC'$), then 
it is not clear whether the  Jacobian syzygy constructed in \eqref{eq12.5} is primitive. Due to this fact, Theorem \ref{thm1} cannot be regarded as a special case of Theorem \ref{thm13}.

\end{rk}

Exactly as in Corollary \ref{corPen1}, when $mdr(f)$ is known, the freeness or nearly freeness of $\CC$ is determined by the global Tjurina number $\tau(\CC)$. This can be seen in  the examples given in the fourth section as well as in the following three results, to be proved in the last two sections.

\begin{thm}
\label{thmTAN}

Let $\HH:h=0$ be an irreducible  curve in $\PP^2$, of degree $e\geq 3$  and having $\delta \geq 0$ nodes and $\kappa\geq 0 $ simple cusps as singularities.
Let $p$ be a generic point in $\PP^2$, such that there are exactly $m_0=e(e-1)-2\delta -3\kappa$ simple tangent lines to $\HH$, say $L_1,...,L_{m_0}$, passing through $p$. 
Assume moreover that the $\delta$  (resp. $\kappa$) secants $L'_j$ (resp. $L''_k$) determined by  the point $p$ and the nodes (resp. the cusps) of $\HH$ are transversal to $\HH$ at each intersection point $q$, i.e. $(L,\HH)_q=\mult_q\HH$ for $L=L'_j$ or $L=L''_k$.
Then the curve 
$$\CC=\HH \cup ( \cup_{i=1,m_0} L_i) \cup (\cup_j L'_j) \cup (\cup_kL''_k)$$
 is free with exponents $(e,e^2-e-1-\delta -2\kappa)$ and the complement $U=\PP^2 \setminus \CC$ is a $K(\pi,1)$-space.

\end{thm}

In the case of line arrangements we have the following result, saying that any line arrangement is a subarrangement of a free, $K(\pi,1)$ line arrangement.

\begin{thm}
\label{thmTAN2}

For any line arrangement $\A$ in $\PP^2$ and any point $p$ of $\PP^2$ not in $\A$, the line arrangement $\B(\A,p)$ obtained from $\A$ by adding all the lines determined by the point $p$ and by each of the multiple points in $\A$ is a free, $K(\pi,1)$ line arrangement.

\end{thm}
The case when $p$ and several multiple points of $\A$ are collinear is allowed, and the line added in such a case is counted just once, hence $\B(\A,p)$ is a reduced line arrangement.
The main part of the next result was  stated and proved by a different method in \cite{JV}, see also the Erratum to that paper.

\begin{thm}
\label{thmPEN}

 Let $\CC:f=0$ be a curve arrangement in $\PP^2$ such that the defining equation has the form
$$f=q_1q_2 \cdots q_m,$$
for some $m \geq 3$, where $\deg q_1=\cdots=\deg q_m=k\geq 2$ and the curves $\CC_i:q_i=0$ for $i=1,...,m$ are members of the pencil spanned by $\CC_1$ and $\CC_2$.
Assume that the pencil $u\CC_1+v\CC_2$ is generic, i.e. the curves $\CC_1$ and $\CC_2$ meet transversely in exactly $k^2$ points.
 Then the following properties are equivalent.

\begin{enumerate}

\item Any singularity of any singular member $C_j^s$ of the pencil $u\CC_1+v\CC_2$ is weighted homogeneous and all these singular members $\CC^s_j$  are among the $m$ curves $\CC_i$ in the curve arrangement $\CC$;

\item The curve  $\CC$ is  free  with exponents $(2k-2, mk-2k+1)$.

\end{enumerate}

When the curve $\CC$ is free, then the complement $U=\PP^2 \setminus \CC$ is a $K(\pi,1)$-space.

\end{thm}

I would like to thank Aldo Conca and Jean Vall\`es for useful discussions related to this paper.

\section{Multiple points and Jacobian syzygies } 

\subsection{Proof of Theorem \ref{thmF}} \label{ss0}

If $\CC$ is free with exponents $(r_0,d-r_0-1)$, then the formula for $\tau(\CC)$ is well known, see for instance \cite{DS14}.

Suppose now that $r_0<r:=mdr(f) \leq (d-1)/2$. Then one has
\begin{equation}
\label{eq55}
(d-1)^2-r_0(d-1-r_0)>\phi_1(r):=(d-1)^2-r(d-1-r) \geq \tau(\CC)
\end{equation}
since the function $\phi_1(r)$ is stricly decreasing on $[0,(d-1)/2]$, and Theorem 3.2 in \cite{duPCTC} yields the last inequality.  Next suppose that $r_0 <r$ and $(d-1)/2<r \leq d-r_0-1$. It follows  from Theorem 3.2 in \cite{duPCTC} that  one has
\begin{equation}
\label{eq53}
\tau(\CC) \leq \phi_2(r):=(d-1)^2-r(d-r-1) -{2r+2-d \choose 2}.
\end{equation}
The function $\phi_2(r)$ is strictly decreasing on $((d-4)/2,+\infty)$  and moreover
$$\phi_1\left(\frac{d-1}{2}\right)=\phi_2\left(\frac{d-1}{2}\right).$$
It follows that in this case one also has
$\tau(\CC) <(d-1)^2-r_0(d-1-r_0).$ Therefore the equality $\tau(\CC) =(d-1)^2-r_0(d-1-r_0)$ holds if and only if $r=r_0$, and one may use \cite{duPCTC} or \cite{Dmax} to
complete the proof of the first claim.

To prove the second claim, consider the closed subvariety $X_r$ in $\PP(S_r^3)\times \PP(S_d)$ given by
$$X_r=\{((a,b,c),f) \  \ : \ \ af_x+bf_y+cf_z=0\}.$$
Note that a polynomial $f \in S_d$ satisfies $mdr(f) \leq r$ if and only if $[f] \in \PP(S_d)$ is in the image $Z_r$ of $X_r$ under the second projection.  If there is $0<r_0\leq (d-1)/2$ such that $\tau=\phi_1(r_0)$, then by the above discussion,  $F(d, \tau)$ is exactly the complement of $Z_{r_0-1} \cap C(d,\tau)$ in $C(d,\tau)$.
If such an $r_0$ does not exist, then $F(d, \tau)=\emptyset$, which completes the proof.

\subsection{Proof of Theorem \ref{thm1}} \label{ss1}

We show first that an intersection point $p$ of multiplicity $m$ gives rise to a syzygy
\begin{equation}
\label{rel_p}
R_p:  a_pf_x+b_pf_y+c_pf_z=0
\end{equation}
where $\deg a_p=\deg b_p=\deg c_p=d-m$ and such that the polynomials $a_p,b_p,c_p$ have no common factor in $S$. Let $f=gh$, where $g$ (resp. $h$) is the product of linear factors in $f$ corresponding to lines in $\A$ passing (resp. not passing) through the point $p$. If we choose the coordinates on $\PP^2$ such that $p=(1:0:0)$, then $g$ is a homogeneous polynomial in $y,z$ of degree $m$, while each linear factor $L$ in $h$ contains the term in $x$ with a non zero coefficient $a_L$. Moreover, $\deg h=d-m$.
It follows that 
$$f_x=gh_x=gh\sum_L\frac{a_L}{L}=f\frac{P}{h},$$
where $P$ is a polynomial of degree $d-m-1$ such that $P$ and $h$ have no common factors.
This implies that
\begin{equation}
\label{syz1}
dhf_x=dPf=xPf_x+yPf_y+zPf_z,
\end{equation}
i.e. we get the required syzygy $R_p$ by setting $a_p=xP-dh$, $b_p=yP$ and $c_p=zP$.

Now, by the definition of $mdr(f)$, we get $mdr(f) \leq d-m$ and it remains to analyse the case
$mdr(f)<d-m.$ Let $R_1$ be the syzygy of degree $mdr(f)$ among $f_x,f_y,f_z$. It follows that
$R_p$ is not a multiple of $R_1$, and hence when
$$\deg R_1 +\deg R_p = mdr(f)+d-m \leq d-1$$
we can use Lemma 1.1 in \cite{ST} and get the case (1). The case (2) is just the situation when the case (1) does not hold, so there is nothing to prove.

\begin{rk}
\label{rklines2}

The method of proof of Theorem \ref{thm1} gives a proof of Theorem \ref{thm13} when $k=1$ and $\HH$  is a reduced curve, not necessarily irreducible. Indeed, $\HH$ reduced implies that $h$ and $h_x$ cannot have any common factor. Any such irreducible common factor would correspond to a line passing through
 $p=(1:0:0)$, and $h$ does not have such factors by assumption.
\end{rk}

 Theorem \ref{thm1} clearly implies the following Corollary, saying that the highest multiplicity of a point of a (nearly) free line arrangement cannot take arbitrary values with respect to the exponents.

\begin{cor}
\label{cor0.5} (i) If $\A$ is a free line arrangement with exponents $d_1 \leq d_2$, then either $m=d_2+1$ or  $m \leq  d_1+1$. 

\noindent (ii) If $\A$ is a nearly free line arrangement with exponents $d_1 \leq d_2$, then either $m=d_2$ or  $m \leq  d_1$.

\end{cor}

The first claim (i) in Corollary \ref{cor0.5} should be compared with the final claim in Corollary 4.5 in  \cite{FV} and looks like a dual result to Corollary 1.2 in \cite{A1}. As a special case of Corollary \ref{cor0.5} we get the following.

\begin{cor}
\label{cor1} (i) If $\A$ is a free line arrangement with exponents $d_1 \leq d_2$ and $m > d/2$, then $m=d_2+1$.

\noindent (ii) If $\A$ is a nearly free line arrangement with exponents $d_1 \leq d_2$ and $m \geq d/2$, then   $m = d_2$. 
\end{cor}

The following consequence of Theorem \ref{thm1} is also obvious.

\begin{cor}
\label{cor1.5} If $\A$ is a  line arrangement such that $2m=d$, then either $mdr(f)=m$ and $\A$ is not free, or $mdr(f)=m-1$ and $\A$ is free with exponents $m-1,m$.
\end{cor}

\subsection{Proof of Proposition  \ref{propDSa}} \label{Saito}

For the reader's convenience, we recall some facts from \cite{DSa}, see also \cite{DS14}.
Let $C$ be a reduced plane curve in $\PP^2$ defined by $f=0$.
Let $\alpha_C$ be the minimum of the Arnold exponents  $\alpha_p$ of the singular points $p$ of $C$. The plane curve singularity $(C,p)$ is weighted homogeneous of type $(w_1,w_2;1)$ with $0<w_j \leq 1/2$, if there are local analytic coordinates $y_1,y_2$ centered at $p$ and a polynomial $g(y_1,y_2)= \sum_{u,v} c_{u,v}y_1^uy_2^v$, with $c_{u,v} \in \C$ and where the sum is over all pairs $(u,v) \in \N^2$ with $uw_1+vw_2=1.$ In this case one has 
\begin{equation}
\label{Ae}
\alpha_p={w_1}+{w_2},
\end{equation}
see for instance \cite{DSa}. 
With this notation, Corollary 5.5 in \cite{DSa} can be restated as follows.

\begin{thm}
\label{vanishing}
Let $C:f=0$ be a degree $d$ reduced curve in $\PP^2$ having only weighted homogeneous singularities. 
Then $AR(f)_r=0$ for all $r <\alpha_C d-2$.
\end{thm}

In the case of a line arrangement $C= \A$, a point $p$ of multiplicity $k$ has by the above discussion the Arnold exponent $\alpha_p=2/k$. It follows that, for $m=m(\A)$, one has
\begin{equation}
\label{AC}
\alpha_C=\frac{2}{m},
\end{equation}
and hence Theorem \ref{vanishing} implies
\begin{equation}
\label{mdr1}
mdr(f) \geq \frac{2}{m}d-2.
\end{equation}
In other words
\begin{equation}
\label{mdr2}
m\geq \frac{2d}{mdr(f)+2},
\end{equation}
i.e. exactly what is claimed in Proposition  \ref{propDSa}.

\subsection{Proof of Corollary \ref{corNF}}
Let $\B$ be defined by $g=0$. Then Corollary \ref{cor1} (ii) applied to $\A$ implies that $d_1=d-m$, and hence in particular
$$\tau(\A)=(d-1)^2-(d-m)(m-1)-1,$$
see \cite{Dmax}.  Note that $\tau(\A)=\tau(\B)$ as for a line arrangement the total Tjurina number is determined by the intersection lattice, see for instance \cite{DStexpo}, section (2.2).
If $mdr(g)=d-m$, then our characterisation of nearly free arrangements in \cite{Dmax} via the total Tjurina number implies that $\B$ is also nearly free.

On the other hand, if $mdr(g)<d-m$, Theorem \ref{thm1} applied to the arrangement $\B$ implies that the only possibility given the assumption $m \geq d/2$ is that $\B$ is free with exponents $m-1, d-m$, in particular
$$\tau(\B)=(d-1)^2-(d-m)(m-1).$$
This is a contradiction with the above formula for $\tau(\A)$, so this case is impossible.

\section{On free line arrangements}  

\subsection{Proof of Corollary \ref{cor2}}

The proof of Corollary \ref{cor2} is based on Theorem 3.1 in \cite{FV}, which we recall now in a slightly modified form, see also \cite{DStexpo}, section (2.2).

\begin{thm}
\label{thmFV}
Let $\B$ be an arrangement of $d$ lines in $\PP^2$ and suppose that there are two integers $k \geq 1$ and $\ell \geq 0$ such that $d=2k+\ell+1$ and there is an intersection point in $\B$ of multiplicity $e$ such that 
\begin{equation}
\label{s1}
k \leq e \leq k+\ell+1.
\end{equation}
Then the arrangement $\B$ is free with exponents $(k, k+\ell)$ if and only if the total Tjurina number of $\B$ satisfies the equality
\begin{equation}
\label{s2}
\tau(\B) =(d-1)^2-k (k+\ell).
\end{equation}
\end{thm}

\begin{rk}
\label{rknewproof}
A new proof of Theorem \ref{thmFV} can be given using our Theorem \ref{thmF}, Theorem \ref{thm1} and  Theorem 3.2 in \cite{duPCTC}.
\end{rk}

To prove the first claim of Corollary \ref{cor2}, we  apply Theorem 2 in \cite{FV} to the arrangement $\B$.
Corollary \ref{cor1} implies that $m=m(\A)=m(\B) \leq d_2+1$. Then we can set $k=d_1$, $
\ell=d_2-d_1$ and $e=m$ and we get 
$$k=d_1 \leq e=m \leq k+\ell +1=d_2+1.$$
It follows that  $\tau(\B) =\tau(\A)=(d-1)^2-d_1d_2$ and hence the line arrangement $\B$ is free with exponents $d_1, d_2$.

\medskip

The last claim of Corollary \ref{cor2} follows since $m<d_1$ implies via Proposition  \ref{propDSa} that 
$$ \frac{2d}{d_1+2} <d_1.$$
But this quadratic inequality in $d_1$ holds if and only if $d_1 > \sqrt {2d+1} -1$.

\subsection{Some examples of  free line arrangements }

Now we consider some examples of  line arrangements.
First we show by examples that all the cases listed in Theorem \ref{thm1} and Corollary \ref{cor1.5} can actually occur inside the class of free line arrangements.

\begin{ex}
\label{ex1}  The line arrangement $$\A: f=xyz(x-z)(x+z)(x-y)=0$$ is free with exponents  $ (2, 3)$ and has $m=4>d/2=3$. Hence we are in the situation $d_1=mdr(f)=2=d-m$.

\end{ex}

\begin{ex}
\label{ex2}  The line arrangement $$\A: f=xyz(x-z)(x+z)(x-y)(x+y)(y-z)=0$$  is free with exponents  $ (3,4)$ and has $m=4=d/2$. Hence we are in the situation $d_1=mdr(f)<4=d-m$ and $d_1=mdr(f)=m-1$, as in Corollary \ref{cor1.5}.

Similarly, the line arrangement $$\A: f=xyz(x-z)(x+z)(x-y)(x+y)(y-z)(y+z)=0$$ is free with exponents  $ (3,5)$ and has $m=4<d/2$. Hence we are in the situation $d_1=mdr(f)<5=d-m$ and $d_1=mdr(f)=m-1$.

\end{ex}

\begin{ex}
\label{ex3}  The line arrangement 
$$\A: f=xyz(x-z)(x+z)(x-y)(x+y)(y-z)(y+z)
(x+2y)(x-2y)$$
$$
(x+2z)(x-2z)(y-2z)(y+2z)
(x+y-z)(x-y+z)(-x+y+z)(x+y+z)=0$$ 
is free with exponents  $ (9,9)$ and has $m=6<19/2$. Hence we are in the situation $m=6\leq mrd(f)=d_1=9 <d-m=13$.

\end{ex}

Finally we give an example showing that the inequality in Proposition  \ref{propDSa} is sharp.

\begin{ex}
\label{ex5} The line arrangement $$\A: f=(x^3-y^3)(y^3-z^3)(x^3-z^3)=0$$ is free with exponents  $ (4,4)$ and has 
$$m=3=  \frac{2d}{d_1+2}.$$

\end{ex}

\section{Pencils and Jacobian syzygies} 

 Let $\CC:f=0$ be a curve arrangement in $\PP^2$ such that the defining equation has the form
$$f=q_1q_2 \cdots q_m h=gh,$$
for some $m \geq 2$, where $\deg q_1=\cdots=\deg q_m=k$ and the curves $\CC_i:q_i=0$ for $i=1,...,m$ are  reduced members of the pencil $\PPP: u\CC_1+v\CC_2$. Assume that $\PPP$ has a $0$-dimensional base locus and that it contains only reduced curves. In terms of equations, one can write 
\begin{equation}
\label{eq11}
q_i=q_1+t_iq_2,
\end{equation}
for $i=3,...,m$ and some $t_i \in \C^*$ distinct complex numbers. In other words, the curve subarrangement $\CC':g=0$ of $\CC$ consists of $m \geq 2$ reduced members of a pencil.

To find a Jacobian syzygy for $f$ as in \eqref{rel_m} is equivalent to finding a homogeneous $2$-differential form $\omega$ on $\C^3$ with polynomial coefficients  such that
\begin{equation}
\label{eq12}
\omega \wedge df=0.
\end{equation}

\subsection{The case when $\CC$ is a pencil }

When $h=1$, i.e. when $\CC=\CC'$ is a pencil itself, then one can clearly take 
\begin{equation}
\label{eq12.5}
\omega =dq_1 \wedge dq_2,
\end{equation}
see also Lemma 2.1 in \cite{JV}. This form yields a primitive syzygy of degree $2k-2$ if we show that

\noindent (i) $\omega \ne 0$, and

\noindent (ii) $\omega$ is primitive, i.e. $\omega$ cannot be written as $e \eta$, for $e \in S$ with $\deg e >0$ and $\eta $ a
$2$-differential form on $\C^3$ with polynomial coefficients. Such a polynomial $e$ is called a {\it divisor} of $\omega$.

The first claim follows from Lemma 3.3 in \cite{Dmax}, since $q_1=0$ is a reduced curve and $q_1$ and $q_2$ are not proportional. The claim (ii) is a consequence of  Lemma 2.5 in \cite{JV}, or can easily be proven directly by the interested reader.

\subsection{The case when $\CC$ is a not pencil }

If $\deg h =d-km>0$, then we set
\begin{equation}
\label{eq13}
\omega=adq_1 \wedge dq_2+bdq_1\wedge dh+cdq_2\wedge dh,
\end{equation}
with $a,b,c \in S$ to be determined.  The condition \eqref{eq12} becomes
$$g\left[ a-bh \left(\frac{1}{q_2}+ \frac{t_3}{q_3}+ \cdots +\frac{t_m}{q_m}\right)+ch\left(\frac{1}{q_1}+ \frac{1}{q_3}+ \cdots +\frac{1}{q_m}\right)\right]dq_1 \wedge dq_2 \wedge dh=0.$$
We have the following result.

\begin{lem}
\label{lem2}
Assume that the curves $\CC_1:q_1=0$, $\CC_2:q_2=0$ and $\HH: h=0$ have no common point. Then
the $2$-form $\omega=adq_1 \wedge dq_2+bdq_1\wedge dh+cdq_2\wedge dh$
with $a=-mh$, $b=-q_2$ and $c=q_1$ is non-zero and satisfies $\omega \wedge df=0$.
Moreover, any divisor of $\omega$ is a divisor of the Jacobian determinant $J(q_1,q_2,h)$ of the polynomials $q_1,q_2,h$ and of $h$. In particular, if $h$ is irreducuble, then $\omega $ is primitive.
\end{lem}

\proof Since the ideal $(q_1,q_2,h)$ is $\bf m$-primary, where ${\bf m}=(x,y,z)$, it follows that
$$dq_1 \wedge dq_2 \wedge dh=J(q_1,q_2,h)dx \wedge dy \wedge dz \ne 0,$$
see \cite{GH}, p. 665. This shows in particular that $\omega \ne 0.$ Indeed, one has 
$$\omega \wedge dq_1= q_1J(q_1,q_2,h)dx \wedge dy \wedge dz$$
 and 
$$\omega \wedge dq_2=q_2J(q_1,q_2,h)dx \wedge dy \wedge dz.$$
Since $q_1$ and $q_2$ have no common factor, these equalities  show that any divisor $e$ of $\omega$ divides $J(q_1,q_2,h)$.

Let $\Delta$ be the contraction of differential forms with the Euler vector field, see Chapter 6 in 
\cite{D1} for more details if needed. Then one has
$$J(q_1,q_2,h)\Delta(dx \wedge dy \wedge dz)=\Delta(dq_1 \wedge dq_2 \wedge dh)=$$
$$=kq_1dq_2 \wedge dh-kq_2dq_1 \wedge dh+(d-mk)hdq_1 \wedge dq_2=$$
$$=k \omega +d \cdot h dq_1 \wedge dq_2.$$
This implies that any divisor $e$ of $\omega$ and of $J(q_1,q_2,h)$ divides $h$ as well.
Since $h$ does not divide $J(q_1,q_2,h)$, see \cite{GH}, p. 659, the last claim follows.
\endproof

The main results based on the above facts are Theorems \ref{thm11} and \ref{thm13}, stated in the Introduction. Their  proofs are exactly the same as the proof of Theorem \ref{thm1} using the discussion above. 

\begin{rk}
\label{rkDISCARD}

 The case $r=d-m$ in Theorem \ref{thm1},  the case $r=2k-2$ in Theorem \ref{thm11}, or the case $r=2k-2+ \deg (h)$ in Theorem \ref{thm13} can  sometimes be discarded if $r=mdr(f) >(d-1)/2$ using the inequality \eqref{eq53}.

\end{rk}

Now we illustrate these results by some examples.

\begin{ex}
\label{ex12} (i) The line arrangement $$\A: f=(x^k-y^k)(y^k-z^k)(x^k-z^k)=0$$ for $k\geq 2$ is seen  to be free with exponents  $ (k+1,2k-2)$ using  Theorem \ref{thmPEN}. This arrangement has $m(\A)=k$ for $k \geq 3$, hence the Jacobian syzygy constructed  in the proof of Theorem \ref{thm1} has degree $d-m(\A)=2k$.
The Jacobian syzygy constructed in \eqref{eq12.5} has degree $d_2=2k-2$, hence we are in the case (1) of Theorem \ref{thm11} when $k \geq 4$. Theorems   \ref{thmF} and \ref{thm1}  give an alternative proof for the freeness of $\A$. The same method shows that the arrangement 
$$\A': f=xyz(x^k-y^k)(y^k-z^k)(x^k-z^k)=0$$ 
for $k\geq 2$ is free with exponents $(k+1,2k+1)$.

(ii) The curve arrangement
$$\CC: f=xyz(x^3+y^3+z^3)[(x^3+y^3+z^3)^3-27x^3y^3z^3]=0$$
is just the Hesse arrangement from \cite{JV} with one more smooth member of the pencil added.
One has $k=3$ and $m=5$, hence $r=mdr(f)=4$ follows from
Theorem \ref{thm11}. Moreover, the Jacobian syzygy constructed in \eqref{eq12.5} has minimal degree $r=mdr(f)$, and this is always the case by Theorem \ref{thm11} when $k=3$ or when $m>3$. To compute the total Tjurina number $\tau(\CC)$ of $\CC$, note that the $9$ base points of the pencil are ordinary $5$-fold points, hence each contributes with $16$ to $\tau(\CC)$. There are four singular members of the pencil in $\CC$, each a triangle, hence we should add $12$ to $\tau(\CC)$ for these $12$ nodes which are the vertices of the four triangles. It follows that
$$\tau(\CC)=9 \times 16+12=156=(d-1)^2-r(d-r-1)=14^2-4\times 10,$$
which shows that $\CC$ is free with exponents $(4,10)$ using \cite{duPCTC}, \cite{Dmax}. 

\end{ex}

\begin{ex}
\label{ex14} 
 (i) The line arrangement $$\A: f=(x^k-y^k)(y^k-z^k)(x^k-z^k)x=0$$ is seen by a direct computation to be free with exponents  $ (k+1,2k-1)$. This arrangement has $m(\A)=k+1$ for $k \geq 2$, hence the Jacobian syzygy constructed  in the proof of Theorem \ref{thm1} has degree $d-m(\A)=2k$.
The Jacobian syzygy constructed in \eqref{eq12.5} has degree $d_2=2k-1$, hence we are in the case (1) of Theorem \ref{thm13}.

(ii) Consider the curve arrangement $\CC:f=x(x^{m-1}-y^{m-1})(xy+z^2)$ for $m \geq 3$.
Here $k=1$ and $d=m+2$. Theorem \ref{thm13} implies that $r=mdr(f)=\deg(h)=2$. 
To compute the total Tjurina number $\tau(\CC)$ of $\CC$, note that $(0:0:1)$ is an ordinary $m$-fold point, hence it contributes to $\tau(\CC)$ by $(m-1)^2$. Each of the $(m-1)$ lines in $x^{m-1}-y^{m-1}$ meets the smooth conic $\HH: xy+z^2=0$ in two points, so has a contribution to 
$\tau(\CC)$ equal to $2$. The line $x=0$ is tangent to $\HH$ at the point $p=(0:1:0)$ and hence at $p$ the curve $\CC$ has an $A_3$ singularity. It follows that
$$\tau(\CC)=(m-1)^2+2(m-1)+3=m^2+2=(d-1)^2-r(d-r-1)-1.$$
Using \cite{Dmax}, we infer that the curve $\CC$ is nearly free with exponents $(2,m)$.
The same method shows that the curve arrangement $\CC':f=xy(x^{m-2}-y^{m-2})(xy+z^2)$ for $m \geq 3$ is free with exponents $(2,m-1)$.

\end{ex}

\subsection{Proof of Theorem \ref{thmTAN}}

 It is known that the degree of the dual curve $\HH^*$ is given by $m_0=e(e-1)-2\delta-3\kappa$, see \cite{GH}, p. 280, hence the existence of points $p$ as claimed is clear.
We apply Theorem \ref{thm13} to the curve $\CC$, with $k=1$, $d=m+e$, where $m=m_0+\delta +\kappa$ is the total number of lines in $\CC$.
Using the known inequality for an irreducible curve
$$ \delta +\kappa \leq \frac{(e-1)(e-2)}{2},$$
one obtains $r=mdr(f)=e$. 

To  compute of the global Tjurina number $\tau(\CC)$ as in Example \ref{ex14} (ii), one has to add the following contributions.

\begin{enumerate}

\item  $\tau(\CC,p)=(m-1)^2=(e^2-e-1-\delta -2\kappa)^2$;

\item   the singularities of $\CC$ along each tangent line $L_i$ except $p$ have a total Tjurina number
$e+1$, so we get in all a contribution 
$$m_0(e+1)= (e(e-1)-2\delta-3\kappa)(e+1)$$
from the tangent lines.

\item  the singularities of $\CC$ along each secant  line $L_j'$ except $p$ have a total Tjurina number
$e+2$, so we get in all a contribution 
$\delta(e+2)$.

\item  the singularities of $\CC$ along each secant  line $L''_k$ except $p$ have a total Tjurina number
$e+3$, see for instance \cite[Lemma 2.7]{STo}, so we get in all a contribution 
$\kappa(e+3)$. Indeed, all the singularities with Milnor number at most 7 are known to be weighted homogeneous, and hence the Tjurina number coincides with the Milnor number in these cases.

\end{enumerate}
When we add up these contributions, we get
$$\tau(\CC)=(d-1)^2-r(d-r-1).$$
Hence $\CC$ is free with exponents $(r,d-r-1)=(e,e^2-e-1-\delta-2\kappa)$ using \cite{duPCTC}, \cite{Dmax}. The fact that the complement $U$ is aspherical follows from the fact that the projection from $p$ on a generic line $L$ in $\PP^2$ induces a locally trivial fibration
$F \to U \to B$, where both the fiber $F$ and the base $B$ are obtained from $\PP^1$ by deleting finitely many points, see \cite[Chapter 4]{D1}.

\subsection{Proof of Theorem \ref{thmTAN2}}

Let $e$ be the number of lines in $\A$ and $m$ the number of extra lines contained in $\B(\A,p)$. We set $d=e+m$.
First we determine $\tau(\B(\A,p))$. Since all the singularities of a line arrangement are weighted homogeneous, then $\tau(\B(\A,p))=\mu(\B(\A,p))$, the total Milnor number of the line arrangement $\B(\A,p)$. This number enters in the following formula for the Euler number of (the curve given by the union of all the lines in ) $\B(\A,p)$.
$$\chi(\B(\A,p))=\chi(C_d)+\mu(\B(\A,p))= 3d-d^2+\mu(\B(\A,p)),$$
where $C_d$ denotes a smooth curve of degree $d$, see \cite[Corollary 5.4.4]{D1}. By the additivity of the Euler numbers, we get
$$\chi(U)=\chi(\PP^2)-\chi(\B(\A,p))=3- 3d+d^2-\mu(\B(\A,p)),$$
where $U=\PP^2 \setminus \B(\A,p)$. Projection from $p$ induces a locally trivial fibration, with total space $U$, and fiber (resp. base) a line $\PP^1$ with $e+1$ (resp. $m$) points deleted.
The multiplicative prperty of the Euler numbers yield
$$\chi (U)=(1-e)(2-m).$$
The last two equalities can be used to get the following
$$\tau(\B(\A,p))=\mu(\B(\A,p))=(d-1)^2-e(m-1).$$
By Theorem \ref{thmF}, it remains only to show that $mdr(f) \geq \min \{e,m-1\}.$ Note that 
Theorem \ref{thm1} applied to the arrangement $\B(\A,p)$ and the multiple point $p$ shows that either $mdr(f)=d-m=e+m-m=e$, or  $mdr(f) \geq m-1$. The proof is complete, since the
 asphericity of $U$ follows exactly as in the proof above, using the obvious locally trivial fibration given by the central projection from $p$.

Note that in this proof the point $p$ is not necessarily the point of highest multiplicity of the arrangement $\B(\A,p)$.

\section{The case of generic pencils} 

Let $\CC:f=0$ be a curve arrangement in $\PP^2$ such that the defining equation has the form
$$f=q_1q_2 \cdots q_m,$$
for some $m \geq 2$, where $\deg q_1=\cdots=\deg q_m=k$ and the curves $\CC_i:q_i=0$ for $i=1,...,m$ are  members of the pencil  $\PPP:u\CC_1+v\CC_2$. 
We say that the pencil $\PPP$ is generic if the following  condition is satisfied: the curves $\CC_1$ and $\CC_2$ meet transversely in exactly $k^2$ points. If this holds, then the generic member of $\PPP$ is smooth, and any member of the pencil $\PPP$ is smooth at any of the $k^2$ base points.
Let us denote by $\CC^s_j$ for $j=1,...,p$ all the singular members in this pencil $\PPP$.
One has the following result.

\begin{prop}
\label{propPEN} If the pencil $\PPP$ is generic, then the sum of the total Milnor numbers of the singular members $\CC^s_j$ in the pencil satisfies
$$\sum_{j=1,p} \mu (\CC^s_j)=3(k-1)^2.$$
 
\end{prop}

\proof First recall that $\mu (\CC^s_j)$ is the sum of the Milnor numbers of all the singularities of the curve $\CC^s_j$. Then we consider two smooth members $D_1: g_1'=0$ and $D_2: g_2'=0$ in the pencil and consider the rational map 
$\phi: X \to \C,$
where $X=\PP^2 \setminus D_1$ and $$\phi(x:y:z)=\frac{g_2'(x,y,z)}{g_1'(x,y,z)}.$$
Then it follows that $\phi$ is a tame regular function, see \cite{SiTi}, whose singular points are exactly the union of the singular points of the curves $\CC^s_j$ for $j=1,...,p$.
From the general properties of tame functions it follows that
$$\sum_{j=1,p} \mu (\CC^s_j)=\sum_{a \in X}\mu(\phi,a)= \chi(X,X\cap D_2 ).$$
Since the Euler characteristic of complex constructible sets is additive we get
$$\chi(X,X\cap D_2 )=\chi(\PP^2)-\chi(D_1) -\chi(D_2)+\chi(D_1 \cap D_2)=$$
$$=3+2k(k-3)+k^2=3(k-1)^2.$$

\endproof

\begin{rk}
\label{rkPEN}
Let $k\geq 2$ and consider the discriminant hypersurface $\D_k \subset \PP(S_k)$ consisting of singular plane curves of degree $k$. Then it is known that $\deg \D_k=3(k-1)^2$, see for instance \cite{JV}. It follows that a generic pencil $\PPP$ as defined above, and thought of as a line in $\PP(S_d)$  has the following {\it transversality property}.

$(T)$ For any intersection point $p \in \PPP \cap \D_k$, one has an equality
\begin{equation}
\label{eq23}
\mult_p\D_k = (\D_k, \PPP)_p,
\end{equation}
where $\mult_p\D_k$ denotes the multiplicity of the hypersurface $\D_k$ at the point $p$ and $(\D_k, \PPP)_p$ denotes the intersection multiplicity of  the hypersurface $\D_k$ and the line $\PPP$ at the point $p$. 

\medskip

To see this, we use the inequality
$\mult_p\D_k \leq (\D_k, \PPP)_p,$
which holds in general and the equality
$\mult_p\D_k =\mu(\CC(p)),$
where $\CC(p)$ is the degree $k$ reduced curve corresponding to the point $p$, see \cite{DD}.
Then one has
$$3(k-1)^2=\deg \D_k=\sum_{p\in \D_k \cap \PPP}(\D_k, \PPP)_p\geq \sum_{p\in \D_k\cap\PPP}\mult_p\D_k=\sum_{p\in \D_k \cap \PPP}\mu(\CC(p))=3(k-1)^2,$$
where the last equality follows from Proposition \ref{propPEN}.
\end{rk}

With the notation above, we also have the following.

\begin{cor}
\label{corDISC}
If $\PPP$ is a generic pencil of degree $k$ plane curves, with $k \geq 2$, then the number of points in the intersection $\PPP\cap \D_k$ satisfies
$$|\PPP\cap \D_k|\geq 3.$$
Moreover the equality $|\PPP\cap \D_k|= 3$ holds if and only if each of the three singular fibers of the pencil $\PPP$ is a union of $k$ concurrent lines, i.e. we are essentially in the situation of Example \ref{ex12} (i).
\end{cor}

\proof
This claim  follows from Proposition \ref{propPEN} and the well known inequality
$$\mu(\CC') \leq (k-1)^2,$$
for any reduced plane curve $\CC'$ of degree $k$, where equality holds if and only if $\CC'$ is a 
union of $k$ concurrent lines. Indeed, for the inequality one can use the primitive embedding of lattices
given in the formula (4.1), p. 161 in \cite{D1} or \cite[Proposition 7.13]{Loo}. When we have equality, the claim follows from the fact that a Milnor lattice of an isolated hypersurface singularity cannot be written as an orthogonal direct sum of sublattices, see \cite[Proposition 7.5]{Loo} for the precise statement.
\endproof
Note that in any pencil $\PPP$ of degree $k$ curves, the number of completely reducible fibers $\CC'$ (i.e. fibers $\CC'$ such that $\CC'_{red}$ is a  line arrangement) is at most $4$ and the only known example is the Hesse pencil generated by a smooth plane cubic and its Hessian, see \cite{Yu}.

\subsection{Proof of Theorem \ref{thmPEN}}

 First we assume (1) and prove (2).  For this, we compute the total Tjurina number $\tau(\CC)$, taking into account the fact that the singularities of $\CC$ are of two types: the ones coming from the singularities of the singular members $\CC^s_j$ and the $k^2$ base points, each of which is an ordinary $m$-fold point. It follows that
\begin{equation}
\label{eq54}
\tau(\CC) =\sum_{j=1,p} \tau (\CC^s_j)+k^2(m-1)^2=3(k-1)^2+k^2(m-1)^2,
\end{equation}
since $\tau (\CC^s_j)=\mu (\CC^s_j)$, all the singularities being weighted homogeneous.

Assume first that $m \geq 4$. Then Corollary \ref{corPen1} implies that $r=mdr(f)=2k-2$
and the equation \eqref{eq54} yields $\tau(\CC)=(d-1)^2-r(d-1-r)$, i.e. $\CC$ is free.

Consider now the case $m=3$. If $r=mdr(f)=2k-2$, the same proof as above works.
Moreover, if we are in the case (1) of Theorem \ref{thm11}, i.e. $mdr(f)=k+1=r_0$, then again we get
$$\tau(\CC)=(d-1)^2-r_0(d-1-r_0),$$
and hence $\CC$ is free in this case as well. It remains to discuss the case (2) in Theorem \ref{thm11}.
 This can be done using Theorem \ref{thmF}, thus completing the proof of the implication $(1) \Rightarrow (2)$. The implication $(2) \Rightarrow (1)$ is obvious using \cite{KS0}.

To prove the last claim, note that the pencil $u\CC_1+v\CC_2$ induces a locally trivial fibration $F \to U \to B$ where the fiber $F$ is a smooth plane curve minus $k^2$ points, and  the base is obtained from $\PP^1$ by deleting finitely many points.

\end{document}